\documentclass[letterpaper, 10 pt, conference]{ieeeconf}  

\IEEEoverridecommandlockouts                              
\usepackage[utf8]{inputenc}
\usepackage{xcolor}
\usepackage{xspace}
\usepackage{textcomp}

\usepackage{multirow}

\usepackage{amsmath,amssymb}
\usepackage{bm,fixmath} 

\usepackage{graphicx}

\usepackage{accents}

\usepackage{flushend}

\usepackage{algorithm}
\usepackage{algpseudocode}
\algrenewcommand\algorithmicrequire{\textbf{Input:}}
\newcommand{\CommentState}[1]{\Statex\hspace{\algorithmicindent}{\color{blue}// #1}}

\DeclareMathOperator*{\argmin}{arg\,min}

\DeclareMathOperator{\proj}{proj}

\newtheorem{remark}{Remark}

\newtheorem{proposition}{Proposition}

\newtheorem{assumption}{Assumption}
\newtheorem{lemma}{Lemma}
\newtheorem{corollary}{Corollary}

\newcommand{\ubar}[1]{\underaccent{\bar}{#1}}

\newcommand{\norm}[1]{\left\lVert#1\right\rVert}
\newcommand{\ev}[1]{\mathbb{E}\left[#1\right]}
\newcommand{\round}[1]{\left\lfloor#1\right\rceil}

\newcommand{\N}{\mathbb{N}}
\newcommand{\R}{\mathbb{R}}

\newcommand{\C}{\mathcal{C}}

\newcommand{\e}{\mathbold{e}}
\newcommand{\x}{\mathbold{x}}
\newcommand{\y}{\mathbold{y}}
\newcommand{\w}{\mathbold{w}}
\newcommand{\z}{\mathbold{z}}

\newcommand{\1}{\pmb{1}}

\def\lmin{{\ubar{\lambda}}}
\def\lmax{{\bar{\lambda}}}

\title{\LARGE \bf
Asynchronous Distributed Learning\\with Quantized Finite-Time Coordination}

\author{Nicola Bastianello, Apostolos I. Rikos, Karl H. Johansson 
\thanks{This work was partially supported by the European Union’s Horizon Research and Innovation Actions programme under grant agreement No. 101070162, and partially by Swedish Research Council Distinguished Professor Grant 2017-01078 Knut and Alice Wallenberg Foundation Wallenberg Scholar Grant.}
\thanks{N. Bastianello and K. H. Johansson are with the School of Electrical Engineering and Computer Science, and Digital Futures, KTH Royal Institute of Technology, Sweden, {\tt\small \{nicolba | kallej\}@kth.se}.}%
\thanks{Apostolos~I.~Rikos is with the Artificial Intelligence Thrust of the Information Hub, The Hong Kong University of Science and Technology (Guangzhou), Guangzhou, China. 
He is also affiliated with the Department of Computer Science and Engineering, The Hong Kong University of Science and Technology, Clear Water Bay, Hong Kong, China. 
E-mail: {\tt~apostolosr@hkust-gz.edu.cn}.}
}

\begin{document}

\maketitle

\thispagestyle{empty}
\pagestyle{empty}


\begin{abstract}
In this paper we address distributed learning problems over peer-to-peer networks. In particular, we focus on the challenges of quantized communications, asynchrony, and stochastic gradients that arise in this set-up.
We first discuss how to turn the presence of quantized communications into an advantage, by resorting to a \textit{finite-time, quantized coordination} scheme. This scheme is combined with a distributed gradient descent method to derive the proposed algorithm.
Secondly, we show how this algorithm can be adapted to allow asynchronous operations of the agents, as well as the use of stochastic gradients.
Finally, we propose a variant of the algorithm which employs zooming-in quantization.
We analyze the convergence of the proposed methods and compare them to state-of-the-art alternatives.
\end{abstract}

\section{Introduction}\label{sec:intro}
In recent years, multi-agent systems have become ubiquitous in a broad range of applications, \textit{e.g.} robotics, power grids, traffic networks \cite{molzahn_survey_2017,nedic_distributed_2018}.
A multi-agent system consists of autonomous agents with communication and computation capabilities, cooperating to accomplish a specific goal, \textit{e.g.} learning, decision-making, navigation. 
In this paper, we will focus on developing algorithms to enable decentralized learning.
In decentralized learning, the agents in the system collect and locally store data, with the goal being to collectively train a model without sharing these raw data \cite{alghunaim2022unified}.
To enable this objective, the design of distributed learning (or optimization) algorithms has been extensively studied in the past decades \cite{nedic_distributed_2018,yang_survey_2019}. In particular, different classes of algorithms have been proposed, with the main ones being gradient methods (\textit{e.g.} DGD), gradient tracking, and dual methods (\textit{e.g.} ADMM) \cite{notarstefano_distributed_2019}.
In this paper, we will focus on gradient-based methods.

Learning over a multi-agent system, however, presents a number of practical challenges, with communication constraints being a central one.
These constraints may arise due to the reliance on communication channels with limited bandwidth (\textit{e.g.} wireless) \cite{qian_distributed_2022}, or the necessity to share high dimensional models (\textit{e.g.} neural networks) \cite{richtarik_3pc_2022}.
A common solution to reduce the communication burden is the use of \textit{quantization}, which however may result in lower accuracy of the trained model.
In this paper we address distributed learning with quantized communication, and aim at showing how to turn quantization from a design constraint into an opportunity.
The central idea is that agents employing quantized communications can \textit{reach an inexact consensus in finite time}.
Thus, in this paper we combine a Finite-Time, Quantized Coordination (FTQC) scheme with gradient descent, to design efficient learning algorithms that only require quantized communications. 
In particular, the FTQC scheme we analyze is based on the consensus ADMM \cite{bastianello_asynchronous_2021}, differently from previous alternatives \cite{2023:Rikos_Johan_IFAC,bastianello_online_2023}.

Besides limited communications, in this paper we address two additional challenges that arise in distributed learning: \textit{asynchrony} \cite{assran_advances_2020} and \textit{stochastic gradients} \cite{xin_decentralized_2020}.
Indeed, the cooperating agents may have \textit{access to different, and limited, hardware resources}. On the one hand, different resources result in the agents having different computation speeds, which make asynchronous completion of local training steps inevitable. In this paper, therefore, we follow the literature \cite{zhao_asynchronous_2015,xu_convergence_2018,tian_achieving_2020} in designing a gradient-based learning algorithm that enables asynchronous local training.
On the other hand, limited hardware resources have the consequence that agents may find the computation of local gradients prohibitive (e.g., due to potentially lengthy computation times or memory constraints). 
For this reason, the agents may resort to computing inexact stochastic gradients \cite{xin_decentralized_2020,yi_primaldual_2022}, by only using a subset of the available data. In the following we design an algorithm that relies on stochastic gradients.

The main contributions of the paper are as follows: 
\begin{enumerate}
    \item We propose an asynchronous distributed learning algorithm which relies on finite-time, quantized coordination. The novel FTQC scheme we propose is based on consensus ADMM, and the algorithm allows for the use of stochastic gradients.

    \item We further propose an alternative version of the algorithm which employs \textit{zooming-in quantization}, which progressively reduces the loss of accuracy due to quantized communications.

    \item We analyze the convergence of the proposed FTQC scheme, and of the complete algorithm, highlighting the effect of (i) quantization, (ii) stochastic gradients, (iii) asynchrony. We further analyze the effect of zooming-in quantization on the convergence.
    
    \item We conclude with numerical results comparing the proposed FTQC scheme and algorithms against state-of-the-art alternatives.
\end{enumerate}

\section{Problem Formulation}\label{sec:prob_form}
Given the undirected graph $\mathcal{G} = (\mathcal{V}, \mathcal{E})$ with $N$ agents, our goal is to solve the \textit{distributed optimization problem}
\begin{equation}\label{eq:problem-consensus}
	\min_{\x \in \R^{n N}} \sum_{i = 1}^N f_i(x_i) \quad \text{s.t.} \ \x \in \C,
\end{equation}
where $\x = [x_1^\top, \ldots, x_N^\top]^\top$, $f_i : \R^n \to \R$ are local costs, each privately held by one of the agents, and we define the \textit{consensus set}
$
    \C = \{ \x \in \R^{n N} \ | \ x_i = x_j \ \forall i, j \in \mathcal{V} \}.
$
In the following we are interested in the \textit{finite-sum local costs} that arise in learning applications, hence we assume that
\begin{equation}\label{eq:local-costs}
    f_i(x) = \frac{1}{m_i} \sum_{h = 1}^{m_i} \ell(x; d_i^h)
\end{equation}
where $\ell : \R^n \to \R$ is a loss function (\textit{e.g.} logistic) and $\{ d_i^h \}_{h = 1}^{m_i}$ are the data points stored by agent $i$ (\textit{e.g.} pairs of label and feature vector).

The following assumptions will hold throughout the paper.

\begin{assumption}[Network]\label{as:network}
The graph $\mathcal{G} = (\mathcal{V}, \mathcal{E})$ is undirected and connected.
\end{assumption}

\begin{assumption}[Costs]\label{as:costs}
The local cost $f_i : \R^n \to \R$ is $\lmin$-strongly convex and $\lmax$-smooth for each agent $i \in \mathcal{V}$.
\end{assumption}

By Assumption~\ref{as:network}, the graph is connected, which ensures that problem~\eqref{eq:problem-consensus} can be solved in a distributed fashion.
Moreover, Assumption~\ref{as:costs} implies that there is a unique solution to the problem, which we can write as
$$
    \x^* = \1_N \otimes x^*, \quad x^* = \argmin_{x \in \R^n} \sum_{i = 1}^N f_i(x).
$$

\smallskip

We are now ready to discuss the objectives that will guide the algorithm design in Section~\ref{sec:algorithm}:
\begin{enumerate}
    \item \textit{Quantized communications}: learning problems are \textit{high-dimensional}, as the model being trained may have a large number of parameters $n \gg 1$ \cite{li_federated_2020,gafni_federated_2022}. However, distributed learning requires the agents to share their local models, which may cause a large communication overhead. 
    The idea is to design an algorithm that uses quantized/compressed communications \cite{zhao_towards_2023}.

    \item \textit{Stochastic gradients}: in order to train an accurate model, the local costs~\eqref{eq:local-costs} of the learning problem are often defined over a large data-set, with $m_i \gg 1$ \cite{li_federated_2020,gafni_federated_2022}. 
    However, computing the gradients of such costs may be excessively time consuming.
    Hence, we are interested in designing an algorithm that uses less computationally expensive gradients, called stochastic gradients.

    \item \textit{Asynchrony}: synchronizing all agents in the network $\mathcal{G}$ may not be feasible, especially when $N \gg 1$ \cite{assran_advances_2020}. 
    The goal is to design an algorithm that allows the agents to perform computations asynchronously.
\end{enumerate}

\section{Algorithm Design}\label{sec:algorithm}
In this section we design the proposed distributed learning algorithm tailored to the objectives detailed in section~\ref{sec:prob_form}.
Conceptually, one could think of solving problem~\eqref{eq:problem-consensus} by applying the \textit{projected gradient} method \cite{taylor_exact_2018} characterized by
\begin{equation}\label{eq:algorithm}
	\x_{k+1} = \proj_\C\left( \x_k - \alpha \nabla f(\x_k) \right), \quad k \in \N,
\end{equation}
where
$
	\nabla f(\x) = [\nabla f_1(x_1)^\top, \ldots, \nabla f_N(x_N)^\top]^\top
$
collects the local gradients, and the projection onto the consensus space is
$
	\proj_\C(\x) = \frac{1}{N} \sum_{i = 1}^N x_i.
$

Clearly, the computation of $\proj_\C(\x)$ cannot be performed in a distributed fashion, except with specific architectures such as federated learning \cite{gafni_federated_2022}. The objective therefore is to propose a distributed (and approximate) implementation of the consensus projection.
Different techniques have been explored to this end, foremost of which is averaged consensus. In particular, we can replace $\proj_\C(\x)$ with one or more consensus steps, giving rise to Near-DGD \cite{berahas_balancing_2019}
\begin{equation}\label{eq:near-dgd}
    \x_{k+1} = \mathbold{W}^t \left( \x_k - \alpha \nabla f(\x_k) \right), \quad k \in \N, \ t \geq 1,
\end{equation}
where $\mathbold{W}$ is a symmetric, doubly stochastic matrix. 
Alternatively, the average consensus can be replaced with dynamic average consensus, which gives rise to gradient tracking algorithms \cite{xin_decentralized_2020}.

\smallskip

In this paper we take a different approach by using, similarly to \cite{2023:Rikos_Johan_IFAC,bastianello_online_2023}, a \textit{finite-time, quantized coordination} (FTQC) scheme to approximate $\proj_\C(\x)$.
Indeed, as discussed in section~\ref{sec:prob_form}, in learning applications we may need to use quantized communications, and the idea is to use this fact to our advantage.

\subsection{Finite-time, quantized coordination}\label{subsec:ftc}
The main insight guiding our design is that \textit{specific consensus schemes achieve convergence in finite-time when the communications are quantized}. Employing such a scheme therefore allows the agents to approximate $\proj_\C(\x)$ in a finite number of iterations.
The algorithm proposed in  \cite{2022:Rikos_Johan_Split_Autom}, for example, is specifically tailored to achieve this goal.
However, we explore a different strategy by showing how the consensus ADMM \cite{bastianello_asynchronous_2021} satisfies the requirements of a FTQC scheme.

Let $\{ y_i \}_{i \in \mathcal{V}}$ be local states that need to be averaged. We can formulate this as the distributed optimization problem
\begin{equation}
	\min_{\x \in \R^{n N}} \frac{1}{2} \sum_{i = 1}^N \norm{x_i - y_i}^2 \quad \text{s.t.} \ \x \in \C,
\end{equation}
to which we apply the distributed ADMM \cite{bastianello_asynchronous_2021}, yielding Algorithm~\ref{alg:finite-time-consensus}
\footnote{To be precise, Algorithm~\ref{alg:finite-time-consensus} is derived from \cite{bastianello_asynchronous_2021} by setting $\alpha = 0.5$, and excluding the termination step, discussed in the following.}.

\begin{algorithm}[!ht]
\caption{Finite-time quantized coordination (FTQC)}
\label{alg:finite-time-consensus}
\begin{algorithmic}[1]
	\Require The states to be averaged $\{ y_i \}_{i \in \mathcal{V}}$, $z_{ij}^0 = 0$ for all $i \in \mathcal{V}$, $j \in \mathcal{N}_i$, penalty $\rho > 0$, quantizer $q(\cdot)$, termination threshold $\theta > 0$.
    \Statex{\color{blue}// initialization}
    \State each agent $i \in \mathcal{V}$ picks $z_{ij}^0$ for all $j \in \mathcal{N}_i$
	\For{$\ell = 1, 2, \ldots$}
    \CommentState{local update and transmission}
    \If{agent $i$ is active}
        \State computes
            $
                w_i^\ell = \frac{1}{1 + \rho |\mathcal{N}_i|} \left( y_i + \sum_{j \in \mathcal{N}_i} z_{ij}^\ell \right)
            $
        \State and transmits $t_{i \to j} = q\left( -z_{ij}^\ell + 2\rho w_i^\ell \right)$ to each neighbor $j \in \mathcal{N}_i$
    \EndIf

    \CommentState{auxiliary update}
    \If{agent $i$ is active and receives $t_{j \to i}$}
        \State computes $z_{ij}^{\ell+1} = \frac{1}{2} \left( z_{ij}^\ell + t_{j \to i} \right)$
    \EndIf

    \CommentState{termination}
    \If{$\norm{z_{ij}^{\ell+1} - z_{ij}^\ell} \leq \theta$ for all $j \in \mathcal{N}_i$}
        \State agent $i$ terminates
    \EndIf
    \EndFor
\end{algorithmic}
\end{algorithm}

The following lemma shows how Algorithm~\ref{alg:finite-time-consensus} can indeed serve as a FTQC scheme. The proof is reported in Appendix~\ref{app:proof-lemma-ftqc}.

\begin{lemma}[Consensus ADMM as FTQC scheme]\label{lem:ftqc-convergence}
Let $\{ w_i^\ell \}_{\ell \in \N}$ be the trajectory generated by Algorithm~\ref{alg:finite-time-consensus} applied to average $\{ y_i \}_{i \in \mathcal{V}}$,  with a given (auxiliary) initial condition $\{ z_{ij}^0 \}_{i \in \mathcal{V}, \ j \in \mathcal{N}_i}$, and penalty $\rho > 0$.
Assume that communications are quantized according to
\begin{equation}\label{eq:symmetric-quantizer}
    q(x) = \Delta \round{\frac{x}{\Delta}}, \ \Delta > 0
\end{equation}
where $\round{\cdot}$ rounds to the nearest integer.
Then there exist $\mu \in (0, 1)$ and $C > 0$ such that for each $i \in \mathcal{V}$
\begin{equation*}
    \norm{w_i^\ell - \frac{1}{N} \sum_{i \in \mathcal{V}} y_i} \leq C \left( \mu^\ell d(z_0) + \frac{\Delta}{2} \sqrt{n \sum_{i \in \mathcal{V}} |\mathcal{N}_i|} \frac{1 - \mu^\ell}{1 - \mu} \right)
\end{equation*}
with $d(z_0)$ being a function of the initial conditions.
Additionally, convergence is achieved after a finite number of iterations; that is, for $\ell \geq \bar{\ell}$ we have $w_i^\ell = w_i^{\ell+1}$, with
$$
    \bar{\ell} \geq \left| \frac{1}{\log(\mu)} \log\left( \frac{\frac{\Delta}{2} \sqrt{n \sum_{i \in \mathcal{V}} |\mathcal{N}_i|}}{(1 - \mu) d(z_0)} \right) \right|.
$$
\smallskip
\end{lemma}

\smallskip

We can now use the finite-time convergence result of Lemma~\ref{lem:ftqc-convergence} to design the termination technique in Algorithm~\ref{alg:finite-time-consensus}. The idea is for each agent $i \in \mathcal{V}$ to detect when the difference $z_{ij}^{\ell + 1} - z_{ij}^\ell$ is below a threshold $\theta$, identifying that their values have stopped changing significantly. In practice, we can choose $\theta = c \Delta$ for some $c \geq 1$.
Notice that the agents do not need to know $\bar{\ell}$ to apply the termination.

\begin{remark}[Speed and accuracy trade-off]\label{rem:admm-ftqc}
Lemma~\ref{lem:ftqc-convergence} shows how the smaller the quantization level $\Delta$ is, the smaller the consensus error.
On the other hand, smaller values of $\Delta$ imply that a larger number of iterations is required to reach convergence, thus presenting a trade-off between speed and accuracy.
\end{remark}

\begin{remark}[Why choose ADMM?]\label{rem:why-admm}
Why choose consensus ADMM as an FTQC scheme, as opposed to the average consensus of Near-DGD, or the FTQC \cite{2022:Rikos_Johan_Split_Autom}? The answer is that ADMM has been proved to be robust to many different challenges, ranging from asynchrony and packet losses \cite{bastianello_asynchronous_2021}, to quantization and other additive errors \cite{bastianello_deplano_online_2023}. Alternative schemes instead lack such theoretical robustness guarantees.
\end{remark}

\begin{remark}[Extensions of Algorithm~\ref{alg:finite-time-consensus}]\label{rem:extensions-ftqc}
Besides allowing for asynchronous activations and packet losses (cf. Remark~\ref{rem:why-admm}), we can further modify Algorithm~\ref{alg:finite-time-consensus} to allow the agents to use different quantizers.
Indeed, Lemma~\ref{lem:ftqc-convergence} would apply equally, but replacing $\Delta$ with the maximum of the local quantization level $\Delta_i$.
\end{remark}

\subsection{Algorithm}\label{subsec:algorithm}
The proposed Algorithm~\ref{alg:main-algorithm} is based on the projected gradient descent~\eqref{eq:algorithm}, where the projection is approximated with the finite-time, quantized coordination scheme discussed in section~\ref{subsec:ftc} above.

In particular, the agents apply a local gradient step in steps 2-3. They then apply Algorithm~\ref{alg:finite-time-consensus} on the result of steps 2-3 ($\y_k$), and update their local states $\x_k$ with the result.
Notice that the algorithm allows for asynchronous operations: steps 2-3 are performed only by active agents, while inactive ones do not update their $y_{i,k}$'s.
Additionally, the agents may use stochastic gradients $\hat{\nabla} f_i$ instead of the full gradients.

\begin{algorithm}[!ht]
\caption{Proposed algorithm}
\label{alg:main-algorithm}
\begin{algorithmic}[1]
	\Require For each agent $i \in \mathcal{V}$ initialize $x_{i,0}$; choose the step-size $\alpha < 2 / \lmax$.
	\For{$k = 0, 1, \ldots$ each agent $i$}
    	\CommentState{local update}
        \If{agent $i$ active}
        \State apply the local (possibly inexact) gradient step
        $$
            y_{i,k} = x_{i,k} - \alpha \hat{\nabla} f_{i}(x_{i,k})
        $$
        \State\hspace{-1.5em}\textbf{else} $y_{i,k} = y_{i,k-1}$
        \EndIf
        \CommentState{coordination}
        \State apply finite-time, quantized coordination
        $$
            \x_{k+1} = \text{Algorithm~\ref{alg:finite-time-consensus}}(\y_k)
        $$
        with $\y_k = [y_{1,k}^\top, \ldots, y_{N,k}^\top]^\top$
    \EndFor
\end{algorithmic}
\end{algorithm}

\subsection{Zooming-in quantization}\label{subsec:zooming-in}
By the discussion in Remark~\ref{rem:admm-ftqc}, the quantization level $\Delta$ mediates the trade-off between the speed of convergence of Algorithm~\ref{alg:finite-time-consensus} and its consensus error.
The idea then is to exploit this trade-off to improve the performance of Algorithm~\ref{alg:main-algorithm} by changing $\Delta$ over time.

By Remark~\ref{rem:extensions-ftqc} we know that the agents can have \textit{uncoordinated} quantizers, \textit{i.e.} $q_i(x) = \Delta_i \round{x / \Delta_i}$.
Each agent then is allowed to modify its quantizer whenever necessary. 
Algorithm~\ref{alg:main-algorithm-zi} reports a prototype of how this can be implemented. Specifically, each agent checks periodically if its local solution $x_{i,k}$ has stopped improving, and selects $r \Delta_i$, $r \in (0, 1)$, if this is the case.

An alternative algorithm with zooming-in quantization was proposed in \cite{rikos_distributed_zooming_2023}. However, in \cite{rikos_distributed_zooming_2023} the agents reduce their quantization in a synchronized fashion via voting, while in Algorithm~\ref{alg:main-algorithm-zi} the agents can set their quantization levels independently.

\begin{algorithm}[!ht]
\caption{Proposed algorithm (zooming-in quantization)}
\label{alg:main-algorithm-zi}
\begin{algorithmic}[1]
	\Require For each agent $i \in \mathcal{V}$ initialize $x_{i,0}$; choose the step-size $\alpha$. Choose the local quantization level $\Delta_i$, and let $r \in (0, 1)$ and $T \geq 1$.
	\For{$k = 0, 1, \ldots$ each agent $i$}
    	\CommentState{local update}
        \If{agent $i$ active}
        \State apply the local (possibly inexact) gradient step
        $$
            y_{i,k} = x_{i,k} - \alpha \hat{\nabla} f_{i}(x_{i,k})
        $$
        \State\hspace{-1.5em}\textbf{else} $y_{i,k} = y_{i,k-1}$
        \EndIf
        \CommentState{coordination}
        \State $
            \x_{k+1} = \text{Algorithm~\ref{alg:finite-time-consensus}}(\y_k),
        $
        with local quantization levels $\Delta_i$
        \CommentState{zooming-in quantization}
        \If{agent $i$ activated $T$ times \& $\norm{x_{i,k+1} - x_{i,k}} \leq \Delta_i$}
            \State $\Delta_i \gets r \Delta_i$
        \EndIf
    \EndFor
\end{algorithmic}
\end{algorithm}

\section{Convergence Analysis}\label{sec:convergence}
In this section we analyze the convergence of Algorithms~\ref{alg:main-algorithm} and~\ref{alg:main-algorithm-zi} when the agents operate asynchronously and apply stochastic gradients.
Before presenting our analysis we make the following assumption.

\begin{assumption}[Set-up]\label{as:random-setup}
Each agent $i \in \mathcal{V}$ activates at iteration $k$ to perform a local gradient step with probability $p_i \in (0, 1]$.
In particular, active agents use a (possibly inexact) gradient $\hat{\nabla} f_i$, for which there exists $\tau \geq 0$ such that
$$
    \ev{\norm{\hat{\nabla} f_i(x) - \nabla f_i(x)}} \leq \tau.
$$
\end{assumption}

\smallskip

Under this assumption, we derive the following convergence result, proved in Appendix~\ref{app:proof-mean-convergence}.

\begin{proposition}[Convergence of Algorithm~\ref{alg:main-algorithm}]\label{pr:mean-convergence}
Let $\{ \x_k \}_{k \in \N}$ be the trajectory generated by Algorithm~\ref{alg:main-algorithm}. Let Assumptions~\ref{as:network}, \ref{as:costs}, and~\ref{as:random-setup} hold.
Then for all $k > 0$ it holds that
\begin{align*}
    \ev{\norm{\x_k - \x^*}} &\leq \sqrt{\frac{\max_i p_i}{\min_i p_i}} \Bigg( \chi^k \norm{\x_0 - \x^*} \\ &+ \left( \gamma + \alpha \tau \sqrt{N} \right) \frac{1 - \chi^k}{1 - \chi} \Bigg)
\end{align*}
where $\chi = \sqrt{1 - (1 - \zeta^2) \min_i p_i} \in (0, 1)$ with $\zeta = \max\{ |1 - \alpha \lmin|, |1 - \alpha \lmax| \}$, and $\gamma = O\left( \Delta \right)$ (as characterized in Lemma~\ref{lem:ftqc-convergence}).
\end{proposition}

\smallskip

As a consequence of Proposition~\ref{pr:mean-convergence}, we see that
$$
    \lim_{k \to \infty} \ev{\norm{\x_k - \x^*}} \leq  \sqrt{\frac{\max_i p_i}{\min_i p_i}} \ \frac{\gamma + \alpha \tau \sqrt{N}}{1 - \chi}
$$
which highlights how the different challenges of quantization, stochastic gradients, and asynchrony impact the asymptotic error.

We can similarly characterize the asymptotic error when we employ the zooming-in quantization of Algorithm~\ref{alg:main-algorithm-zi}. The proof is reported in Appendix~\ref{app:proof-mean-convergence-zi}.

\begin{corollary}[Convergence of Algorithm~\ref{alg:main-algorithm-zi}]\label{pr:mean-convergence-zi}
Let $\{ \x_k \}_{k \in \N}$ be the trajectory generated by Algorithm~\ref{alg:main-algorithm-zi}. Let Assumptions~\ref{as:network}, \ref{as:costs}, and~\ref{as:random-setup} hold.
Then it holds that
\begin{align*}
    \lim_{k \to \infty} \ev{\norm{\x_k - \x^*}} &\leq \sqrt{\frac{\max_i p_i}{\min_i p_i}} \ \frac{\alpha \tau \sqrt{N}}{1 - \chi}.
\end{align*}
\end{corollary}

\smallskip

Clearly, using zooming-in quantization implies that quantization will not impact the asymptotic error, and only the effects of asynchrony and stochastic gradients are present.

\section{Numerical Results}\label{sec:numerical}
In this section we evaluate the performance of the proposed algorithms on a classification task, and compare it with algorithms from the literature.
We consider problem~\eqref{eq:problem-consensus} with local costs
\begin{equation}\label{eq:logistic}
    f_i(x) = \sum_{h = 1}^{m_i} \log\left( 1 + \exp(- b_i^h a_i^h x) \right) + \frac{\epsilon}{2} \norm{x}^2
\end{equation}
defined by the local dataset $\{ d_i^h = (a_i^h, b_i^h) \in \R^{1 \times n} \times \{-1, 1\} \}_{h = 1}^{m_i}$.
In our experiments we have $N = 10$ agents with $m_i = 150$ data-points each, and the problem size is $n = 10$. The regularization weight is set to $\epsilon = 0.075$.
Moreover, unless otherwise stated we use the symmetric quantizer~\eqref{eq:symmetric-quantizer}.
Finally, the data for the problem are randomly generated using the \texttt{make\_classification} utility of \texttt{sklearn} \cite{scikit-learn}, and all algorithms are implemented in \texttt{tvopt} \cite{bastianello_tvopt_2021}.

\subsection{Performance of Finite-Time, Quantized Coordination schemes}\label{subsec:ftqc-evaluation}
We start by comparing the performance of the proposed Finite-Time, Quantized Coordination scheme Algorithm~\ref{alg:finite-time-consensus} with the scheme proposed in \cite{2022:Rikos_Johan_Split_Autom} and employed in \cite{2023:Rikos_Johan_IFAC,bastianello_online_2023}.
In Table~\ref{tab:ftqc} we compare the two FTQC schemes in terms of consensus error and number of iterations, for different quantization levels. The algorithms are applied to average randomly generated vectors in $\R^{10}$, the size of $x$ in~\eqref{eq:logistic}.
\begin{table}[!ht]
\begin{center}
\caption{Consensus error and iteration number for different quantization levels.}
\label{tab:ftqc}
\begin{tabular}{ccccc}
    \hline
    \multirow{2}{*}{$\Delta$} & 
    \multicolumn{2}{c}{\cite{2022:Rikos_Johan_Split_Autom}} & 
    \multicolumn{2}{c}{Algorithm~\ref{alg:finite-time-consensus}}\\ \cline{2-3}\cline{4-5}
               & Cons. err.& Num. iter. & Cons. err. & Num. iter. \\ \hline
    $10^{-8}$  & $5.31 \times 10^{-8}$ & $115$ & $2.85 \times 10^{-8}$ & $105$  \\
    $10^{-7}$  & $5.13 \times 10^{-7}$ & $106$ & $2.88 \times 10^{-7}$ & $95$   \\
    $10^{-6}$  & $5.32 \times 10^{-6}$ & $96$  & $2.89 \times 10^{-6}$ & $86$   \\
    $10^{-5}$  & $5.30 \times 10^{-5}$ & $86$  & $2.86 \times 10^{-5}$ & $80$   \\
    $10^{-4}$  & $5.31 \times 10^{-4}$ & $75$  & $2.89 \times 10^{-4}$ & $66$   \\
    $10^{-3}$  & $5.27 \times 10^{-3}$ & $67$  & $2.87 \times 10^{-3}$ & $57$   \\
    $10^{-2}$  & $5.36 \times 10^{-2}$ & $55$  & $2.91 \times 10^{-2}$ & $48$   \\
    $10^{-1}$  & $5.18 \times 10^{-1}$ & $49$  & $2.87 \times 10^{-1}$ & $43$   \\
    $1$        & $5.17$                & $38$  & $2.89$                & $29$   \\
    \hline
\end{tabular}
\end{center}
\end{table}
We can see that Algorithm~\ref{alg:finite-time-consensus} consistently outperforms \cite{2022:Rikos_Johan_Split_Autom} in terms of consensus error, since it reaches a smaller neighborhood of the consensus, and in terms of the number of iterations it requires.

\smallskip

\begin{figure}[!ht]
    \centering
    \includegraphics[scale=0.5]{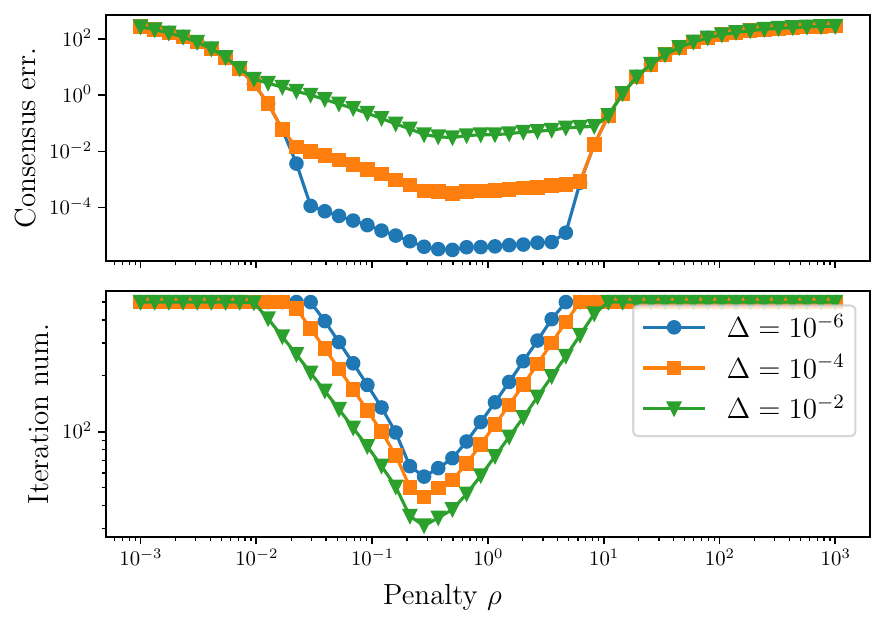}
    \caption{Consensus error and iteration number for Algorithm~\ref{alg:finite-time-consensus} with different penalties and quantization levels.}
    \label{fig:admm-params}
\end{figure}

Turning exclusively to Algorithm~\ref{alg:finite-time-consensus}, we know that it is characterized by two parameters, the penalty $\rho$ and the quantizer $q(\cdot)$. In the following we provide results to guide the tuning of these parameters.
First of all, Figure~\ref{fig:admm-params} reports the consensus error and number of iterations for different values of the penalty and quantization levels.
Interestingly, both metrics are minimized for a value of $\rho \approx 0.3$. Moreover, as predicted by Remark~\ref{rem:admm-ftqc}, the smaller the quantization level is, the smaller the consensus error is, to the detriment of the number of iterations needed to converge.

Finally, Table~\ref{tab:alt-quantization} reports the performance of the consensus scheme with different quantizers besides the symmetric~\eqref{eq:symmetric-quantizer}, namely: floor $q(x) = \Delta \lfloor x / \Delta \rfloor$, ceil $q(x) = \Delta \lceil x / \Delta \rceil$, sparisfier, which sets to zero all components of $x$ with absolute value below $\theta = 0.1$.
\begin{table}[!ht]
\begin{center}
\caption{Performance of Algorithm~\ref{alg:finite-time-consensus} with different quantizers.}
\label{tab:alt-quantization}
\begin{tabular}{ccc}
    \hline
    Quantizer   & Cons. err.                & Num. iter. \\ \hline
    Symmetric   & $3.37 \times 10^{-3}$     & $52$  \\
    Floor       & $8.70 \times 10^{-3}$     & $52$  \\
    Ceil        & $8.69 \times 10^{-3}$     & $52$  \\
    Sparsifier   & $3.79 \times 10^{-3}$     & $50$  \\
    \hline
\end{tabular}
\end{center}
\end{table}
We notice that the floor (employed in \cite{2022:Rikos_Johan_Split_Autom,2023:Rikos_Johan_IFAC,bastianello_online_2023}) and ceiling quantizers attain a more than double the consensus error of the symmetric one. The sparsifier instead achieves similar performance. Future work will explore the use of the sparsifier from a theoretical perspective.

\subsection{Comparison of gradient descent schemes}
The previous section evaluated the performance of the Finite-Time, Quantized Coordination scheme Algorithm~\ref{alg:finite-time-consensus} which is used as a building block of Algorithm~\ref{alg:main-algorithm}.
In this section we discuss the performance of Algorithm~\ref{alg:main-algorithm} itself, and compare it with alternative methods.

\begin{figure}[!ht]
    \centering
    \includegraphics[scale=0.5]{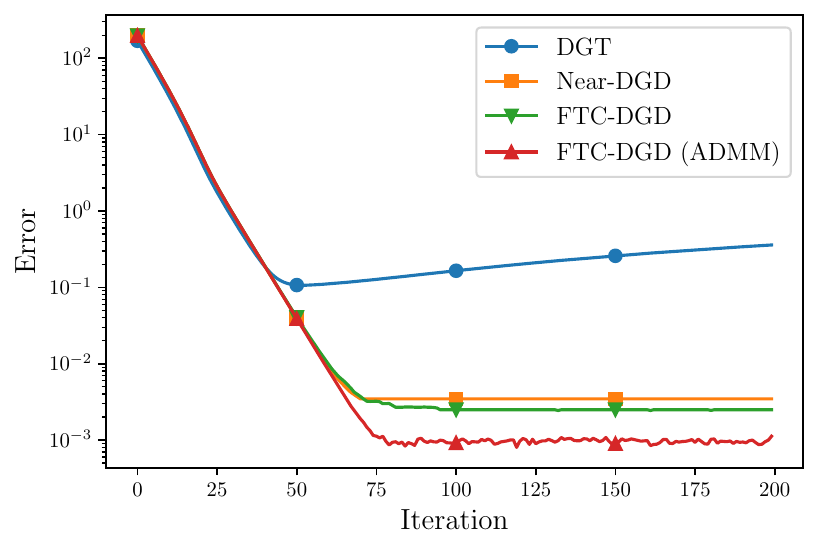}
    \caption{Comparison of different distributed optimization methods with quantized communications.}
    \label{fig:comparison}
\end{figure}

We start by comparing Algorithm~\ref{alg:main-algorithm} to FTQC-DGD \cite{bastianello_online_2023}, Near-DGD \cite{berahas_balancing_2019}, and the distributed gradient tracking (DGT) method \cite{xin_decentralized_2020}. The latter two do not employ a finite-time coordination scheme, but they are modified to use multiple rounds of communications to match the budget of FTQC-DGD and Algorithm~\ref{alg:main-algorithm}.
Figure~\ref{fig:comparison} reports the error trajectories of all methods.
We can see that Algorithm~\ref{alg:main-algorithm} achieves a smaller asymptotic error than Near-DGD and FTQC-DGD, owing to the improved coordination performance of Algorithm~\ref{alg:finite-time-consensus} (cf. section~\ref{subsec:ftqc-evaluation}). Moreover, DGT appears to diverge, which is known to happen with some gradient tracking schemes perturbed by (quantization) noise \cite{bin_stability_2022}.

Table~\ref{tab:comparison} further compares Near-DGD, FTQC-DGD and Algorithm~\ref{alg:main-algorithm} for different quantization levels.
\begin{table}[!ht]
\begin{center}
\caption{Comparison of Near-DGD \cite{berahas_balancing_2019}, FTQC-DGD \cite{bastianello_online_2023}, and Algorithm~\ref{alg:main-algorithm} for different quantization levels.}
\label{tab:comparison}
\begin{tabular}{cccc}
    $\Delta$    & Near-DGD \cite{berahas_balancing_2019}    & FTQC-DGD \cite{bastianello_online_2023}      & Algorithm~\ref{alg:main-algorithm} \\
    \hline
    $10^{-10}$  & $3.75 \times 10^{-8}$     & $4.76 \times 10^{-8}$    & $2.13 \times 10^{-8}$    \\
    $10^{-8}$   & $3.00 \times 10^{-7}$     & $7.35 \times 10^{-7}$    & $1.04 \times 10^{-7}$    \\
    $10^{-6}$   & $3.29 \times 10^{-5}$     & $2.84 \times 10^{-5}$    & $7.79 \times 10^{-6}$    \\
    $10^{-4}$   & $3.46 \times 10^{-3}$     & $2.49 \times 10^{-2}$    & $1.12 \times 10^{-3}$    \\
    $10^{-2}$   & $1.88 \times 10^{-1}$     & $2.56 \times 10^{-1}$    & $7.86 \times 10^{-2}$    \\
    $1$         & $24.02$                   & $20.15$                  & $3.46$                   \\
    \hline
\end{tabular}
\end{center}
\end{table}
The proposed Algorithm~\ref{alg:main-algorithm} outperforms both alternatives, again owing to the improved coordination precision.

\subsection{Variations of Algorithm~\ref{alg:main-algorithm}}
In this section we discuss the performance of Algorithm~\ref{alg:main-algorithm} in challenging scenarios, and compare it to that of Algorithm~\ref{alg:main-algorithm-zi}.

\begin{figure}[!ht]
    \centering
    \includegraphics[scale=0.5]{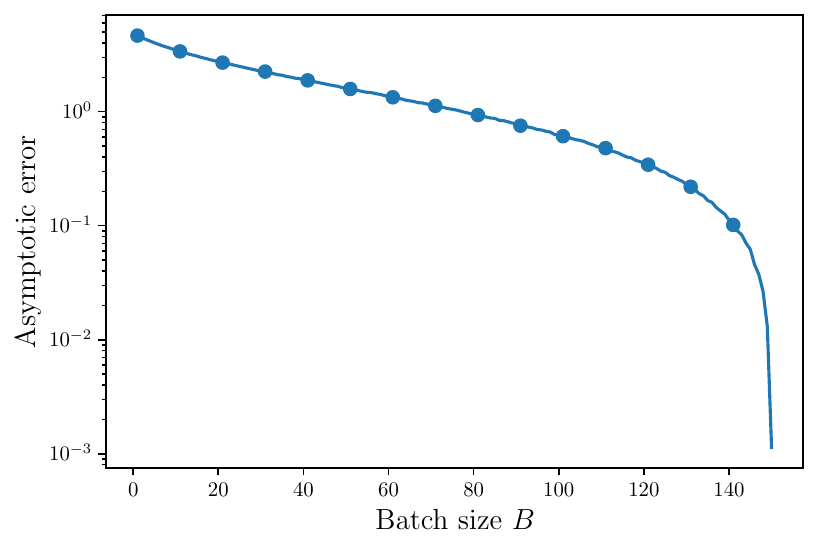}
    \caption{Asymptotic error of Algorithm~\ref{alg:main-algorithm} when the agents employ stochastic gradients of different batch sizes.}
    \label{fig:stochastic-grad}
\end{figure}

As discussed in section~\ref{sec:prob_form}, in learning problems the use of full gradients may be prohibitive, and the agents need to resort to stochastic gradients.
In Figure~\ref{fig:stochastic-grad} we report the asymptotic error achieved by Algorithm~\ref{alg:main-algorithm} when stochastic gradients computed on different batch sizes $B$ are used.
Clearly, the larger the batch size, the better the performance. However, due to the use of quantization, even with $B = m_i$ the algorithm can only reach a neighborhood of the optimal solution.

\begin{figure}[!ht]
    \centering
    \includegraphics[scale=0.5]{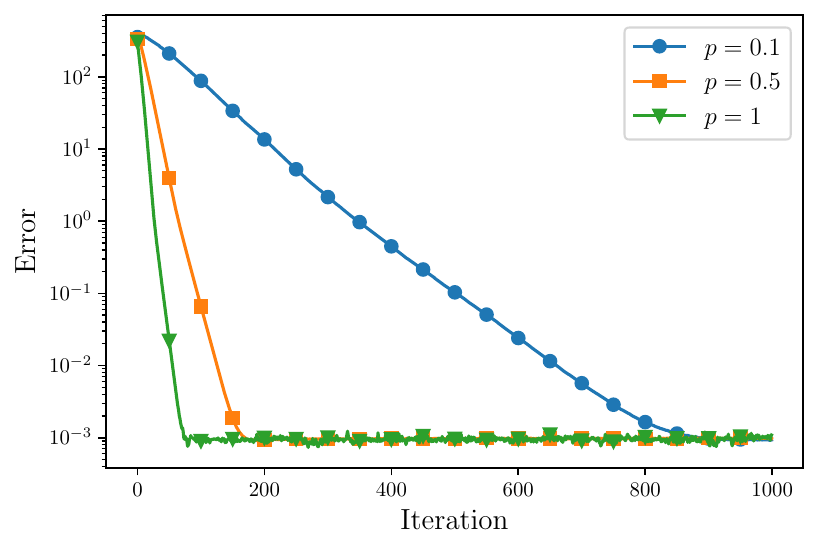}
    \caption{Error trajectory of Algorithm~\ref{alg:main-algorithm} with different agent activation probabilities.}
    \label{fig:asynchronous}
\end{figure}

Another one of the challenges discussed in section~\ref{sec:prob_form} is the asynchronous operation of the agents.
In Figure~\ref{fig:asynchronous} we report the performance of Algorithm~\ref{alg:main-algorithm} in this scenario, when agents activate to perform a gradient descent step with probability $p \in (0, 1]$.
As predicted by the theory \cite{bastianello_stochastic_2023}, the smaller $p$ is the fewer updates are performed, and hence the slower the convergence is.

\begin{figure}[!ht]
    \centering
    \includegraphics[scale=0.5]{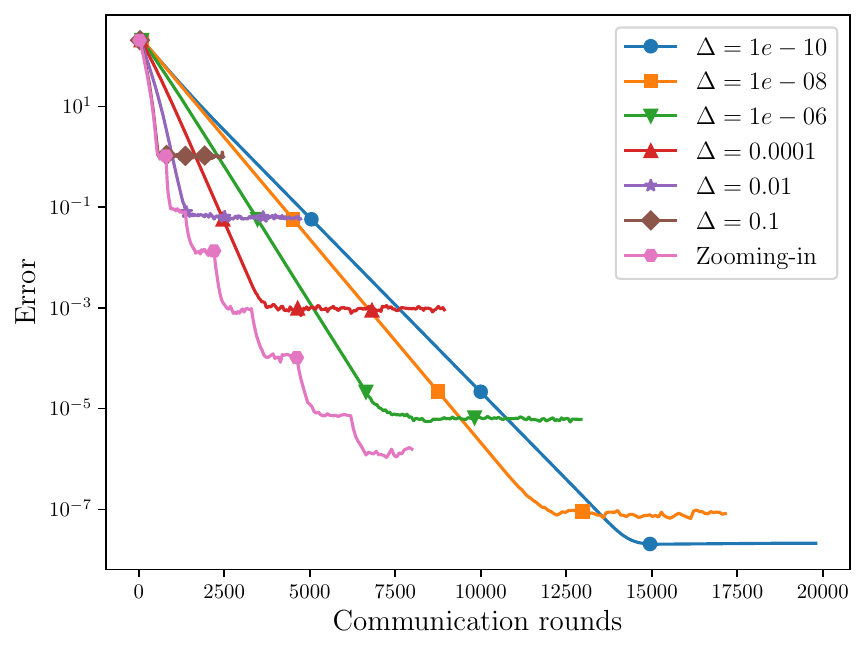}
    \caption{Comparison of Algorithm~\ref{alg:main-algorithm} (fixed quantization level) with Algorithm~\ref{alg:main-algorithm-zi} (zooming-in quantization).}
    \label{fig:zooming-in}
\end{figure}

We conclude this section by comparing the performance of Algorithm~\ref{alg:main-algorithm} with the variation Algorithm~\ref{alg:main-algorithm-zi} that employs zooming-in quantization ($T = 25$, $r = 0.1$).
In particular, Figure~\ref{fig:zooming-in} depicts the error trajectory of the latter against the error trajectory of the former with different quantization levels. The x-axis marks the cumulative number of communication rounds.
We can thus deduce that using zooming-in quantization can achieve very good performance (in terms of asymptotic error) with a smaller number of communication rounds.

\section{Conclusions}\label{sec:conclusions}

In this paper we addressed distributed learning problems over peer-to-peer networks, with a particular focus on the challenges of quantized communications, asynchrony, and stochastic gradients that arise in this set-up.
We first discussed how to turn the presence of quantized communications into an advantage, by resorting to a \textit{finite-time, quantized coordination} scheme. This scheme is combined with a distributed gradient descent method to derive the proposed algorithm.
Secondly, we showed how this algorithm can be adapted to allow asynchronous operations of the agents, as well as the use of stochastic gradients.
Finally, we proposed a variant of the algorithm which employs zooming-in quantization.
We analyzed the convergence of the proposed methods and compared them to state-of-the-art alternatives. The performance of the proposed methods compares very favorably with the alternatives from the literature.



\appendices

\section{Proof of Lemma~\ref{lem:ftqc-convergence}}\label{app:proof-lemma-ftqc}
We start by observing that Algorithm~\ref{alg:finite-time-consensus} consists of an affine update in $\z = [ z_{ij} ]_{i \in \mathcal{V}, \ j \in \mathcal{N}_i}$; in particular, for appropriate matrices and vectors we can write $\z^{\ell+1} = \mathbold{T} \z^\ell + \mathbold{u} + \e^\ell$, $\w^\ell = \mathbold{H} \z^\ell$. The vector $\e^\ell$ represents the noise caused by quantization, that is $e_{ij}^\ell = q(-z_{ji}^\ell + 2 \rho w_j^\ell) - (-z_{ji}^\ell + 2 \rho w_j^\ell)$.
Since Algorithm~\ref{alg:finite-time-consensus} is an affine operator (plus additive noise) then it is $\mu$-metric subregular for a given $\mu \in (0, 1)$ \cite{themelis_supermann_2019}.
Therefore the assumptions of \cite[Theorem~3]{bastianello_deplano_online_2023} are verified, and we have
\begin{align}
    &d(\z^\ell) \leq \mu^\ell d(\z^0) + \sum_{h = 0}^{\ell - 1} \mu^{\ell - h - 1} \norm{\e^h} \label{eq:convergence-1} \\
    &\quad \norm{\w^\ell - \frac{1}{N} \sum_{i \in \mathcal{V}} y_i \otimes \1_N} \leq C d(\z^\ell), \label{eq:convergence-2}
\end{align}
where $d(\z)$ measures the distance of $\z$ from the set of fixed points $\{ \bar{\z} \ | \ \bar{\z} = \mathbold{T} \bar{\z} + \mathbold{u} \}$.

Now, since $\e^\ell$ represents the quantization noise, we can upper bound its norm as follows:
\begin{align*}
    \norm{\e^\ell}^2 &= \sum_{i \in \mathcal{V}} \sum_{j \in \mathcal{N}_i} \norm{q(-z_{ji}^\ell + 2 \rho w_j^\ell) - (-z_{ji}^\ell + 2 \rho w_j^\ell)}^2 \\
    &\leq \sum_{i \in \mathcal{V}} \sum_{j \in \mathcal{N}_i} n (\Delta / 2)^2 = n (\Delta / 2)^2 \sum_{i \in \mathcal{V}} |\mathcal{N}_i|
\end{align*}
where the inequality holds because the quantization commits an error of at most $\Delta / 2$.
Using this bound and combining \eqref{eq:convergence-1} with \eqref{eq:convergence-2} then yields the first thesis.

The goal now is to show that Algorithm~\ref{alg:finite-time-consensus} achieves finite-time convergence.
By \eqref{eq:convergence-1}, we know that $\lim_{\ell \to \infty} d(\z^\ell) = \frac{1}{1 - \mu} \frac{\Delta}{2} \sqrt{n \sum_{i \in \mathcal{V}} |\mathcal{N}_i|}$. Then to bound the time of convergence we impose that the first term on the right-hand side of~\eqref{eq:convergence-1} be smaller than $\lim_{\ell \to \infty} d(\z^\ell)$.
By rearranging, taking the logarithm and the absolute value, the thesis follows. $\hfill\square$

\section{Proof of Proposition~\ref{pr:mean-convergence}}\label{app:proof-mean-convergence}
Algorithm~\ref{alg:main-algorithm} was derived in section~\ref{sec:algorithm} as an inexact version of the projected gradient method, where Algorithm~\ref{alg:finite-time-consensus} replaces the projection onto the consensus set.
Additionally, by Assumption~\ref{as:random-setup}, the agents apply inexact gradients during local computations.
Accounting for both these sources of errors, we can characterize Algorithm~\ref{alg:main-algorithm} as
\begin{equation}\label{eq:algorithm-inexact}
	\x_{k+1} = \proj_\C\left( \x_k - \alpha \nabla f_k(\x_k) \right) + \e_k^q + \e_k^g,
\end{equation}
where $\e_k^q$ is the error due to Algorithm~\ref{alg:finite-time-consensus} (cf. Lemma~\ref{lem:ftqc-convergence}), and $\e_k^g$ is the error due to inexact gradients:
\begin{align*}
    \e_k^p &= \text{Algorithm~\ref{alg:finite-time-consensus}}(\y_k) - \proj_\C(\y_k) \\
    \e_k^g &= \proj_\C(\x_k - \alpha \hat{\nabla} f_k(\x_k)) - \proj_\C(\x_k - \alpha \nabla f_k(\x_k)).
\end{align*}
Moreover, Assumption~\ref{as:random-setup} allows the agents to activate asynchronously, each with its probability $p_i \in (0, 1]$. This means that the $i$-th coordinate of $\x_k$ is updated with probability $p_i$.

Finally, we notice that by the choice $\alpha < 2 / \lmax$, the projected gradient method (without errors and asynchrony) is $\zeta = \max\{ |1 - \alpha \lmin|, |1 - \alpha \lmax| \}$-contractive \cite{taylor_exact_2018}.
This implies that Algorithm~\ref{alg:main-algorithm} can be interpreted as a projected gradient method with bounded additive noise and random coordinate updates. Thus it verifies the assumptions of \cite[Proposition~1]{bastianello_stochastic_2023}, which implies
\begin{align*}
    \ev{\norm{\x_k - \x^*}} &\leq \sqrt{\frac{\max_i p_i}{\min_i p_i}} \Big( \chi^k \norm{\x_0 - \x^*} \\ &+ \sum_{h = 0}^k \chi^{k - h} \norm{\e_h^p + \e_h^g} \Big).
\end{align*}
Now, by Assumption~\ref{as:random-setup} we know that $\ev{\norm{\e_k^g}} \leq \tau \sqrt{N}$, and by Lemma~\ref{lem:ftqc-convergence} we can bound $\norm{\e_k^p} \leq \sqrt{N} C \frac{\Delta}{2} \sqrt{n \sum_{i \in \mathcal{V}} |\mathcal{N}_i|} \frac{1}{1 - \mu} = O(\Delta)$, and the thesis follows. $\hfill\square$

\section{Proof of Corollary~\ref{pr:mean-convergence-zi}}\label{app:proof-mean-convergence-zi}
Following the same derivation as Appendix~\ref{app:proof-mean-convergence} yields
\begin{align*}
    \ev{\norm{\x_k - \x^*}} &\leq \sqrt{\frac{\max_i p_i}{\min_i p_i}} \Big( \chi^k \norm{\x_0 - \x^*} \\ &+ \sum_{h = 0}^k \chi^{k - h} \norm{\e_h^p} + \norm{\e_h^g} \Big).
\end{align*}
By Assumption~\ref{as:random-setup} we know that $\ev{\norm{\e_k^g}} \leq \tau \sqrt{N}$.
On the other hand, by the use of zooming-in quantization, and by Lemma~\ref{lem:ftqc-convergence}, we have
$$
    \norm{\e_k^p} \leq \sqrt{N} C \frac{\Delta_k}{2} \sqrt{n \sum_{i \in \mathcal{V}} |\mathcal{N}_i|} \frac{1}{1 - \mu}
$$
with $\Delta_k = \max_i \Delta_{i,k}$ being the largest quantization level among all agents at time $k$.
We know that $\Delta_k$ is monotonically non-increasing, and that in particular it decreases at finite intervals, when all agents have stopped seeing an improvement in their local solution $x_{i,k}$ (cf. lines 7-8 in Algorithm~\ref{alg:main-algorithm-zi}). Thus $\lim_{k \to \infty} \norm{\e_k^p} = 0$, and by \cite[Lemma~3.1(a)]{sundharram_distributed_2010} $\lim_{k \to \infty} \sum_{h = 0}^k \chi^{k - h} \norm{\e_h^p} = 0$ since $\chi \in (0, 1)$, and the thesis follows. $\hfill\square$


\bibliography{references}
\bibliographystyle{IEEEtran}

\end{document}